\documentclass[12pt]{amsart}
\usepackage{epsfig,amssymb,epic,eepic,color}
\usepackage[all]{xy}
 \usepackage[T1]{fontenc}        
\newtheorem{thm}{Theorem}[section]

\newtheorem{prop}[thm]{Proposition}
\newtheorem{defn}[thm]{Definition}

\newtheorem{cor}[thm]{Corollary}

\newtheorem{remarque}[thm]{Remark}

\newtheorem{rien}[thm]{}
 
\newcommand{\be}{\begin{enumerate}}
\newcommand{\ee}{\end{enumerate}}
\newcommand{\bi}{\begin{itemize}}
\newcommand{\ei}{\end{itemize}}

\def\R{\mathbb{R}}

\def\Z{\mathbb{Z}}

\def\ga{\gamma}    
\def\Ga{\Gamma}

\def\al{\alpha}
\def\be{\beta}

\def\si{\sigma}

\def\ep{\varepsilon}

\def\nd{\noindent}
\def\bull{\hfill$\Box$}
\def\proof{\nd {\bf Proof.\ }}

\numberwithin{equation}{section}

\begin{document}
%\today
\vskip 1cm
\begin{center}
{\sc The Poincar\'e-Lefschetz pairing viewed on Morse complexes

%in case of non-empty boundary
\bigskip

Fran\c cois Laudenbach

}
\end{center}

\title{}
\author{ }
\address{Laboratoire de
math\'ematiques Jean Leray,  UMR 6629 du
CNRS, Facult\'e des Sciences et Techniques,
Universit\'e de Nantes, 2, rue de la
Houssini\`ere, F-44322 Nantes cedex 3,
France.}
\email{francois.laudenbach@univ-nantes.fr}



\thanks{}

\begin{abstract} Given a compact manifold with a non-empty boundary and
 equipped with a generic Morse function (that is, no critical point on the boundary and the restriction to the boundary
 is a Morse function), we already
knew how to construct two Morse complexes, one yielding the absolute homology and the other the relative homology.
In this note, we construct a short exact sequence from both of them and the Morse complex of the boundary. Moreover,
we define a pairing of the relative Morse complex with the absolute Morse complex which induces the intersection 
product in homology, in the form due to S. Lefschetz. This the very first step in an ambitious approach towards 
$A_\infty$-structures buildt from similar data.
\end{abstract}
\maketitle

\thispagestyle{empty}
\vskip 1cm
\section{Introduction}
\medskip

We are given an $n$-dimensional compact manifold $M$ with a non-empty boundary 
$\partial M$ and a Morse function 
$f: M\to \R$ which is generic with respect to the boundary, meaning that $f$ has no critical point on the boundary 
and that the restriction $f_\partial$ of $f$ to $\partial M$ is a Morse function. 
It is well-known that the set of critical points of $f_\partial$  
is divided into two {\it types} $+$ and $-$:
\begin{equation}
crit f_\partial= crit^+ f_\partial \sqcup crit^- f_\partial\,.
\end{equation}
A point $x$ belongs to 
$crit^+ f_\partial$ (resp. $crit^- f_\partial$) if it is a critical point of $f_\partial$ and
 the differential $df(x)$ is positive (resp. negative) on a tangent vector 
at $x$ pointing outwards.\footnote{Here, we choose to introduce notations which are more suggestive than in \cite{lauden}.}

We have introduced in \cite{lauden} the notion of quasi-gradients\footnote{In \cite{lauden}, these vector fields
are named {\it pseudo-gradients} though they vanish at points of $\partial M$ where $df$ does not vanish. So, we prefer to name them {\it quasi-gradients}.}
%of {\it positive} (resp. {\it negative}) type 
{\it positively} (resp. {\it negatively})
{\it adapted to $f$}. Such vector fields, noted respectively noted $X^+$ and $ X^-$, satisfy: 
 \begin{itemize}
  \item[--] $X^+ $ vanishes only in $crit f\cup crit^+f_\partial$ and $\langle df, X^+\rangle >0$ elsewhere;
 \item[--]$X^-$ vanishes only in $crit f\cup crit^-f_\partial$ and $\langle df, X^-\rangle <0$ elsewhere.
 \end{itemize}
  The zeroes of both of them are assumed  hyperbolic, implying the existence of local stable and unstable 
  manifolds. The quasi-gradient $X^+$ (resp. $X^-$) is required to be tangent to the boundary near $crit^+f_\partial$
  (resp. $crit^-f_\partial$).
  Globally, both $X^+$ and $X^-$ are nowhere pointing outwards along $\partial M$. As a consequence,
  their flows are positively complete, and hence, global unstable manifolds exist. 
  By taking inverse images of the local stable manifolds by the
  positive semi-flow, 
  global stable manifolds are also defined  (see Section \ref{globalstable}). 
  
  These  invariant submanifolds are denoted by $W^s(x, X^\pm)$ and $W^u(x, X^\pm)$ respectively 
  when $x$ is a zero of the considered quasi-gradient. If $x\in crit f$, the dimension of $W^s(x, X^+)$ 
  (resp. $W^u(x,X^-)$)
  is equal to the Morse index of $f$ at $x$. If $x\in crit^- f_\partial$,  the dimension of $W^u(x, X^-)$ is equal to the Morse 
  index of $f_\partial$ at $x$; but, if $x\in crit^+f_\partial$, we have
  \begin{equation}\label{dim}
 \dim W^s(x, X^+)= Ind_xf_\partial +1.
 \end{equation}
  
  It makes sense to assume $X^\pm$ Morse-Smale (mutual transversality of stable and unstable manifolds); 
  this property is open and dense. 
  An orientation is chosen on each stable (resp. unstable) manifold arbitrarily when dealing with $X^+$
  (resp. $X^-$). This makes the unstable (resp. stable) manifolds co-oriented
  and allows us to put a {\it sign} on the orbits in $W^s(x, X^+)\cap W^u(y, X^+)$ when the sum of the
  %gradings of $x$ and $y$ 
   codimensions is equal to $n-1$; and similarly for $X^-$.

 Thus, two Morse complexes 
 $C_*(f, X^+)$ and $C_*(f,X^-)$  are built
 whose homologies are respectively isomorphic to $H_*(M,\partial M; \Z)$ 
 and to $H_*(M; \Z)$.\footnote{For defining the differential of 
 these complexes, only the {\it local} stable manifolds are needed.} 
 By abuse of notation, we first neglect to mention the choice of orientations; this will be corrected  in \ref{abuse} 
 for further need.
 For brevity, they are also noted $C_*^+$ and $C_*^-$.

 To be more precise, $C_k^+$ is freely generated by $crit_kf \cup crit_{k-1}^+f_\partial$ (note the shift in the grading
 due to (\ref{dim})) while $C_k^-$  is freely generated by $crit_kf \cup crit_{k}^-f_\partial$.
 The differential $\partial^+:= \partial^{X^+}$ evaluated on a  generator $x\in C_k^+$
 is given by the algebraic counting of orbits of $X^+$ ending at $x$ and starting from generators of $C_{k-1}^+$.
 And similarly for the complex $C_*^-$.
The present note is aimed at proving two results which are stated below.

\begin{thm} \label{exact}
Let $X_\partial$ be a Morse-Smale descending pseudo-gradient of $f_\partial$ on the boundary
$\partial M$ and let $C_*(f_\partial,X_\partial)$ be the associated Morse complex. Then 
for suitable adapted quasi-gradients $X^-$ and $X^+$, there exist a quasi-isomorphic 
extension 
$\widehat C_*(f,X^-)$ 
of the complex $C_*(f,X^-)$ and a short exact sequence of complexes
\begin{equation}
0\longrightarrow C_*(f_\partial,X_\partial)\longrightarrow\widehat C_*(f,X^-)\
\longrightarrow C_*(f,X^+)\longrightarrow 0.
\end{equation}
\end{thm}

%In \cite{kron}, Kronheimer and Mrwoka have considered a less generic setting of functions and gradients since they are assumed both to extend to the double of $M$. This implies some difficulties for constructing the differential of their complexes; namely, there may exist orbits on $\partial M$ connecting critical points of different types $+$ and $-$ (see \cite{kron}\fl{TBC}). But, this price being paid, these complexes enter in a short exact sequence easily.

The second result is stated right below.
I should add that Theorem \ref{pairing}  corrects something which was poorly said at
the end of \cite{lauden}.

\begin{thm} \label{pairing}Here, $M$ is assumed \emph{oriented}. 
For a generic choice of the adapted quasi-gradients $X^+$ and $X^-$,
there is a pairing at the chain level
$$C_k(f, X^+)\otimes C_{n-k}
(f, X^-)\to \Z$$
which induces the intersection pairing in homology
$$\iota: H_*(M,\partial M;\Z)\otimes H_{n-*}(M; \Z)\to \Z
$$
\end{thm}

Intitially, this note was thought of as the beginning of an article on multiplicative structures, namely 
$A_\infty$-algebra structures,  on  Morse complexes \cite{abba}. It appeared that
the pairing $C_*^+\otimes C_{n-*}^-\to \Z$ was not of the same type in nature as the multiplications
of these %above-evoked 
$A_\infty$-structures. Therefore, I decided to separate this piece from \cite{abba}.

\section{A short exact sequence}

We first describe the suitable adapted quasi-gradients $X^+$ and $X^-$ in   Theorem
\ref{exact}.
Let $X_\partial$ be a vector field on $\partial M $
which is a Morse-Smale descending pseudo-gradient of $f_\partial$ and gives rise to the usual Morse complex
of the boundary $C_*(f_\partial,X_\partial)$; its differential is denoted by
$\partial_{\partial M}$.  
 By partition of unity, one easily constructs
a quasi-gradient $X$ of $f$ which extends $X_\partial$.  
This $X$ is tangent to the boundary, and hence
it is not an {\it adapted quasi-gradient}. But it satisfies $X\cdot f<0$ everywhere 
except at the critical points of $f$ and $f_\partial$, where it vanishes with some non-degeneracy condition.
The flow of $X$ is complete,   positively and negatively as well. Therefore, 
one can make $X$ Morse-Smale. % meaning that  the stable and unstable manifolds are mutually transverse.

When $x\in crit_k^- f_\partial$, the unstable $W^u(x,X)$
 coincides
with $W^u(x,X_\partial)\cong \R^k$ and is contained in the boundary. 
The  stable manifold  $W^s(x, X)$ is diffeomorphic to  $\R^{n-k}_{\geq 0}$ and is bounded by 
$W^s(x,X_\partial)$. In the same way, when $y\in crit_k^+f_\partial$, the unstable manifold
$W^u(y,X)$ coincides with the unstable manifold $W^u (y,X_\partial)\cong \R^{n-1-k}$ and is contained
in the boundary. Moreover, the unstable manifold $W^s(y, X)$ is diffeomorphic to $\R^{k+1}_{\geq 0}$
 and is bounded by $W^s(y,X_\partial)$.
 
\begin{remarque} \label{condition*} Since $X$ is tangent to the boundary
there are no connecting orbits of $X$ descending from $x\in crit^-f_\partial$
to $y\in crit f$. Similarly, there are no connecting orbits of $X$ descending from $x\in crit f$
to $y\in crit^+ f_\partial$. %Of course, all other connecting orbits may happen.
\end{remarque}
 We now change  $X$  to  $X^-= X+Y$, which will be  negatively adapted to $f$,
 just  by adding a small vector field Y which satisfies the following condition: % $(*)$:
 \begin{equation}\label{*}
 \left\{ 
\begin{array}{l}
 1)\ Y\ {\rm vanishes\  on\  a\  closed\ neighborhood\  }U\ {\rm of\  }crit^- f_\partial\  {\rm in\ } M;\\
 2)\   Y\  {\rm  points\  inwards\  along\ } 
  \partial M\smallsetminus U\ {\rm  and\ satisfies\  }Y\cdot f\leq 0\ {\rm everywhere}\, ;\\
 3) \ Y\ {\rm vanishes\ away\  from\ a\ neighborhood \ of\ }\partial M\,.
 \end{array}
 \right.
 \end{equation}
 Similarly, $-X$ can be perturbed to $X^+$, which will be  positively adapted to $f$; just take
 $X^+= -X+Z$ where $Z$ is
 a small vector field vanishing on a neighborhood $V$ of $crit^+f_\partial$ in $M$,
   pointing inwards along $\partial M\smallsetminus V$ and satisfying
 $Z\cdot f\geq0$ everywhere. The perturbations $Y$ and $Z$ are small enough so that Remark \ref{condition*}
still applies. So, $X^-$ and $X^+$ will be the desired quasi-gradients of Theorem \ref{exact}.
 
 \begin{prop} \label{no-D}Assume $crit f^+_\partial$ is empty.
 Then the  Morse complex $C_*(f_\partial, X_\partial)$ embeds as a subcomplex
 of $C_*(f,X^-)$. Moreover, one has the following short exact sequence:
 $$
 0\longrightarrow C_*(f_\partial, X_\partial)\mathop{\longrightarrow}\limits^i 
 C_*(f,X^-)\longrightarrow C_*(f,X^+)\longrightarrow 0\, .
 $$
 \end{prop}

 \proof  The embedding $i$ is induced by the inclusion 
 $$crit f _\partial 
 =crit^- f _\partial \hookrightarrow 
 \left(crit f \cup crit^- f _\partial\right).$$ 
 We have to prove that $i$ is a chain morphism. This will follow
 from equalities (1) and (2) below.
 Let $x\in crit_k f_\partial=crit^-_k f_\partial$.
 %Since $X^-$ is close to $X$, 
 By Remark \ref{condition*} applied to $X^-$, for every $y\in crit_{k-1} f$ we have 
  $$(1) \quad \quad  \langle\partial^- x, y\rangle=0.$$

  If $y\in crit^-_{k-1} f_\partial$,
 the  intersection  $W^u(x,X^-)\cap W^s_{loc}(y,X^-) $, which is transverse in $M$, can be pushed by an $f$-preserving isotopy % is isotopic
  to $W^u(x,X_\partial)\cap W^s_{loc}(y,X_\partial) $, which is a transverse 
 intersection in $\partial M$ -- note that %in a level set of $f$ right above $y$, 
 $W^s_{loc}(y,X_\partial)$ is % lies in
  the boundary of $W^s_{loc}(y, X^-)$. %=W^s_{loc}(y,X)$
 Then, the signed number of connecting orbits is  the same
 for both quasi-gradients and we have
 $$(2)\quad\quad \langle\partial^- x, y\rangle= \langle\partial_{\partial M}x, y\rangle.$$ 
 
 For the exactness of the  sequence, observe that the complex $C_*(f,X^+)$ is generated
  by the critical points of $f$. Both vector fields $X^+$ and $X^-$ are approximations
  of the Morse-Smale vector field $X$ (up to sign). Therefore, for every $x\in crit_k f$
and $y\in crit_{k-1}f$, the signed number of connecting orbits is the same when counted
with $X^-$ or $X^+$:
$$\langle\partial^+ x, y\rangle=\langle\partial^- x, y\rangle.
$$
  The quotient kills $crit^- f_\partial$, which generates the image of $C_*(f_\partial, X_\partial)$,
  and also the connecting orbits from $crit f$   to  $crit^-f_\partial$.
   The exactness follows.\bull\\

  \nd {\bf Proof of Theorem \ref{exact}}
  
  It was shown in \cite{lauden} (Lemma 2.4),\footnote{
   After that \cite{lauden} appeared, I was informed that a similar lemma exists in 
   \cite{morse} in a setting where only the Morse inequalities are discussed.}
   that there is  a $C^0$-small deformation, supported in a neighborhood $U$ of $crit^-f_\partial$, 
   of the generic  Morse function 
  $f$ to a new generic  Morse function $f'$ with the following property:
  %such that 
  each $x\in crit_k^-f_\partial$ becomes a critical point of \emph{positive} type and index 
  $k$. The degree of 
    $x$ as generator of  % \in %crit_{k+1}^+ 
    $C_*(f')$ is $k+1$. This is obtained at the cost of 
   a new critical point $ x'\in int\, M$ for $f'$, of index $k$ and close to $x$.  The two critical points $x$ and 
   $x'$ of $f'$ are indeed linked by a unique gradient line; since $x$ belongs to the boundary, this pair is not 
   cancellable but its fusion cancels $x'$ only and changes the type of $x$ from $+$ to $-$. 
   
   Arguing  this way with the function $-f$ 
   leads to the following. There exists a $C^0$-small
   deformation of $f$, supported in a neighborhood  $V$
   of $crit^+ f_\partial$,
   to some generic function $\hat f$ having the following property: %the same restriction to the boundary
   $\hat f_\partial=f_\partial$ and %such that  
   each $x\in crit_k^+ f_\partial$
   becomes a critical point of negative type  and index $k$, that is,  $x\in crit_k^- \hat f_\partial$.
   This is made at the cost of a critical point $\hat x\in int\, M$ for $\hat f$ of index $k+1$
   and close to $x$ and satisfying %the important inequality:
   $$\hat f (\hat x)>\hat f(x)\,.
   $$
   The extension which is mentioned in Theorem \ref{exact} consists of
   adding to $C_*^-$ a pair of new generators $\{x,\hat x\}$ for each $x\in crit^+ f_\partial$.
   More precisely, 
   $$ \widehat C_*(f,X^-):= C_*(\hat f, \widehat X^-)
   $$
   for some quasi-gradient $\widehat X^-$ negatively adapted to $\hat f$. 
    % Without  computing the differential $\widehat\partial^N$ for this extension $\widehat C_*^N$ 
     %we know that  $\widehat\partial^N\circ\widehat\partial^N=0$ since this complex is an
   %$N$-complex for $\hat f$. 
   According to \cite{lauden}, the new complex is quasi-isomophic  to the 
   old one $C_*(f,X^-)$. Since the restriction $f_\partial=\hat f\vert_{\partial M}$
    has no critical point of positive type, 
    Proposition \ref{no-D} applies and
   there is an exact sequence
   $$0\longrightarrow C_*(f_\partial,X_\partial)\longrightarrow 
   C_*(\hat f,\widehat X^-)\ \longrightarrow C_*(\hat f,\widehat X^+)\longrightarrow 0,
  $$
   %Here $\widehat X^-$ and  
   where $\widehat X^+$ denotes a suitable vector field positively %type
    adapted to $\hat f$. 
   In order to 
   identify the quotient in this exact sequence, it is necessary to specify 
  this vector field $\widehat X^+$.
   
   In its  support $V$, the modification from $f$ to $\hat f$
 is modelled similarly to  the birth 
   of a pair of critical points in usual Morse Theory.
   The model produces also 
    a descending quasi-gradient $\widehat X$ of $\hat f$ from the quasi-gradient $X$
   of $f$, which coincides with $X$ out of $V$ and on $\partial M$.
   Then, $-\widehat X$ (which is tangent to the boundary) is changed to  $\widehat X^+$  by adding
    a  vector field  $Z$ which is small
 with respect to $\widehat X$ and satisfies 
  the  condition %$(*)$ after  Remark \ref{condition*}.\\
(\ref{*}) up to sign.\\

    \nd {\sc Claim.} {\it The bijection $j:crit^+ f_\partial
\cup crit f\to crit  \hat f $ which maps $x\in crit^+f_\partial$
to $\hat x\in crit \hat f$ and which is the identity on $crit f\subset
crit \hat f$
induces a chain  isomorphism
$C_*(\hat f,\widehat X^+)\cong C_*(f,X^+)$\,. }\\

  \nd {\sc Proof of the claim.} Say $x\in crit_k^+f_\partial$.
  %By construction, 
 On the one hand,  each %transverse 
 $X$-orbit descending from $x$ to $y\in crit_kf$ % ascending from $y\in crit_kf$ to $x$ 
   gives rise to an %transverse 
   $\widehat X$-orbit from %$y$ to $\hat x$
 $\hat x$ to $y$ 
 and hence, an $\widehat X^+$-orbit from $y$ to $\hat x$. 
   Similarly,
   each $X$-orbit on $\partial M$ descending
  from $x$ to $y\in crit_{k-1}^+f_\partial$ gives rise to  an 
  $\widehat X^+$-orbit from $\hat y$ to $\hat x$. And conversely. Making $j$ an identification,
 this proves the following:
$$\partial^{X^+}\, x=\partial{\widehat X^+}\, \hat x\,.
$$

On the other hand, we have to consider $y\in crit_{k+2} f$ and compute
its two differentials with respect to $X^+$ and $\widehat X^+$ and evaluate 
them at $x$
 (recall that $x$ has degree $k+1$ in $C^+_*$). When $x$ and $y$ have consecutive critical values,
 as a consequence of
   Remark \ref{condition*},
   there are no %$\widehat
    $X^+$-connecting orbits from %$\hat 
    $x$ to $crit_{k+2} f$. %\hat f$. 
  
    But, if their critical values are not consecutive, one could have 
 %But possibly, there are 
 a broken $X$-orbit from $y$ to 
$x$ made of an orbit
from $y$ to $z\in crit_{k+1}^-f_\partial$ and an orbit from $z$ to $x$ on 
$\partial M$. By using the deformation formula $X^+=-X+Z$,
such a broken orbit gives rise to an $X^+$-orbit from $x$ to $y$, and hence to an
$\widehat X^+$-orbit from $\hat x$ to $y$. Then, such connecting orbits may exist. Conversely, by looking at the fusion of the pair $(x,\hat x)$ we get that every 
$\widehat X^+$-orbit from $\hat x$ to $y$ is produced by an $X^+$-orbit  from $x$ to $y$.
Then, via $j$ the following equality holds true:
$$\langle\partial^{X^+}\, y,x\rangle=\langle\partial^{\widehat X^+}\, y, \hat x\rangle\,.
$$
This finishes the proof of  the claim and  Theorem \ref{exact}
  follows.\bull\\

 \section{Global stable manifolds and application to intersection pairing}\label{globalstable}
  
  We now discuss the question of {\it global stable manifolds } 
 for adapted quasi-gradients.
 We only 
  consider   $X^-$ in the definition below;
there is a similar definition for $X^+$. 
If $x\in \partial M$
  is a critical point of negative type, so far we  have only considered 
its local 
  stable manifold $W^s_{loc}(x, X^-)$. If $x$ is of index $k$, it
is a small half-disc $D^{n-k}_-$
  whose  planar boundary lies in a 
  level set of $f$  and spherical boundary lies in $\partial M$.
 Since  the flow of $X^-$, noted %by $\bar X^-$ and by 
  $ X^-_t$ at time $t$, is  positively complete, 
 the following definition makes sense:
  
  \begin{defn}  For $x\in crit f \cup crit^-(f_\partial)$, the global stable
 manifold  of 
  $x$ with respect to $X^-$ is defined as the union
  $$W^s(x,X^-)= \mathop{\bigcup}\limits_{t>0}
 \left( X_t^-\right)^{-1}\left(W^s_{loc}(x, X^-)\right).
  $$
  \end{defn}
  
  Under mild   assumptions, 
it is a (non-proper) submanifold with boundary and its closure
  is a stratified set. The following assumption (Morse-Model-Transversality)
   is made in what follows.\smallskip
   
  \begin{itemize}
  \item[(MMT)] {\it For every $x\in crit f\cup crit^- f_\partial$ and $y\in crit^- f_\partial$,  the 
  neighborhood $U_y$ of $y$ in $\partial M$ where $X^-$ is tangent to the boundary 
  is mapped by the flow transversely to $W^s_{loc}(x, X^-)$.
     } 
  \end{itemize}
  \smallskip
  
 \nd Notice that if $X^-$ is Morse-Smale, the transversality condition  is satisfied along 
  a small neighborhood $U$ of the local unstable manifold $W^u_{loc}(y, X^-)$. 
Then, after some 
  small perturbation of $X^-$ on $U_y\smallsetminus U$ destroying the tangency 
  of $X^-$ to $\partial M$, condition (MMT) is fulfilled for the pair $(y,x)$.
 Thus, condition (MMT) is generic among the 
  negatively adapted vector fields.

  \begin{prop}${}$ \label{frontier} If the negatively quasi-gradient $X^-$ is Morse-Smale and  fulfils condition {\rm(MMT)}
  then the following holds:
  
  \nd  $1)$   The global stable manifold $W^s(x, X^-)$ is a 
  submanifold  with boundary (non-closed in general); its boundary lies in $\partial M$.
  
   \nd $2)$ If $z$ lies in the frontier of $W^s(x, X^-)$ in $M$, then it belongs to the stable manifold of some 
   critical point $y$ in   $crit f\cup crit^- f_\partial$ such that $\dim W^s(y, X^-)<\dim W^s(x,X^-)$.
      \end{prop}
   
   This statement 
   also holds for stable manifolds of critical points in $crit f\cup crit^+ f_\partial$ with respect 
   to positively adapted vector fields.\\

  \proof  1) According to the Implicit Function Theorem, the conclusion  
  is clear  near any point where 
  $X^-$ is transverse to the boundary. Near a point $z$ of $U_y$, it follows 
  from (MMT).\\
  
  \nd 2)  This fact is well known in the case of closed manifolds. It is an easy consequence of the Morse-Smale 
  assumption. The proof is alike if the boundary is non-empty. \hfill\bull\\
  % (see \cite{laudenbach}); 
  %its proof amounts to look at how the gradient lines behave when crossing a Morse model.
  %For the case of manifolds with non-empty boundary, it is alike, except that half Morse models have to be taken into 
  %account. ${}$\hfill\bull\\

 % \begin{center}
 % \begin{figure}
%  \input tangency.pstex_t
%\vspace{.5cm}
  
 % \caption{}
  %\end{figure}
 % \end{center}

% Let  $F^q= F^q \left(W^s(x, X^-)\right)$ denote the stratum of codimension $q$ in the above closure.For instance, $F^1$ is the set of points which belong to a broken orbit ending at $x$ with a single break at some critical point of index $k+1$. 
  
 \begin{remarque} \label{fundrem}
  {\rm Due to the transversality assumptions, a small 
 perturbation of %$f, 
 $X^-$ (resp. $X^+$) moves %the union
 each of stable and unstable
manifolds %of each of them 
by a small isotopy,
 and hence, preserves the complex $C^-_*$ (resp.  $C^+_*$ ) up to a canonical 
 isomorphism.}
 
 As a consequence, without  changing the above-mentioned complexes, we are allowed to assume
  that the $X^-$-unstable manifolds of
 $crit f\cup crit^-f_\partial$ intersect the global $X^+$-stable manifolds 
 of $crit f\cup crit^+f_\partial$ transversely.
 \\
  \end{remarque}

   \begin{rien} {\bf Where an abusive notation is corrected.}\label{abuse} {\rm If the orientation of some of 
 the unstable manifolds is changed then the differential
 of the considered Morse complex (absolute or relative) 
 is changed by a non-trivial isomorphism. So, to understand the role of
 the orientability of $M$ in what follows, it will be better to replace $C_*(f, X^-)$ with $C_*(f,X^-, \ep_f^-)$
 %at the place where it is crucial.
 where $\ep_f^-$ denotes the chosen {\it orientation map} which associates an orientation 
  of $W^u(x,X^-)$ with each $x\in crit f\cup crit^- f_\partial$. Note that $\ep_f^-$ orients the unstable manifolds regardless
 of  the quasi-gradient since they all have isotopic  germs at the critical points. And similarly
 for $C_*(f, X^+)$. Actually, we will only apply this change of notation  at the places where it will be crucial.\\

 }
 \end{rien}

\begin{rien} {\bf The  Poincar\'e-Lefschetz isomorphism.}
{\rm At the homology level,  this % Poincar\'e duality
isomorphism is a %\emph{natural}
 isomorphism
$$P: H_{*}(M,\partial M;\Z)\to H^{n-*}(M;\Z),$$ 
We wish to describe it %(and prove it) 
by means of our Morse complexes in order to deduce %have later
a Morse theoretical description of the homological intersection. There are several steps to 
achieve.
}
\end{rien}

\nd 1) First, we recall that there is a \emph{natural} isomorphism at the homology level 
$$I_*(f,X^-): H_*\left(C_*(f,X^-)\right)\to H_*(M;\Z).$$
Indeed, we have described in \cite{lauden} a canonical process for removing the 
 critical points of $f_\partial$ of negative
type.  Once this is done, the unstable manifolds of $X^-$
emerging  from $crit(f)$ yield %up to \emph{simple homotopy equivalence}, 
a cell decomposition of $M$ (see \cite{laudenbach}\footnote{ In this reference, a stronger assumption is made 
on the vector field which implies this cell decomposition to be a $CW$-complex. Without this assumption, the cell decompositon has only the homotopy type of a $CW$-complex. This is sufficient for our discussion.})
whose homology is canonically isomorphic to the singular homology of $M$ (see \cite{milnor}, p. 90).

We now explain the naturality of this isomorphism. Let $(g, Y^-)$ be another pair of generic Morse 
function and negatively adapted quasi-gradient. The choice of a generic path $\ga$
from $(g, Y^-)$ to $(f,X^-)$ gives rise to some \emph{simple homotopy equivalence}
\begin{equation}
\ga_*:C_*(g,Y^-)\to C_*(f,X^-) %\,.
\end{equation}
well defined up to the orientations.\footnote{The creation/cancellation times of pair of critical points along $\ga$ do not allow us to carry orientations along the path.}
At each time that $\ga$ crosses a  stratum corresponding to a codimension-one defect of 
genericity of the pair ({\it function, negatively adapted quasi-gradient}) this yields an \emph{elementary}
modification of the Morse complex, indeed a \emph{quasi-isomorphism}  \cite{lauden}. One checks at each 
occurrence that this quasi-isomorphism is compatible to the isomorphism 
with $H_*(M;\Z)$. %And a change of orientation of some unstable manifolds 
Finally, $\ga_*$ is the composition of all these quasi-isomorphisms. It induces  an isomorphism $[\ga_*]$ in homology making the next diagram commute: 
\begin{equation}\label{natural}
\end{equation}
\vskip -1.3cm{
{\[\xymatrix@1{
 H_*\left(C_*(g,Y^-)\right) \ar[r]^(0.5){[\ga_*]}\ar[rd]_{I_*(g,Y^-)}&
H_*\left(C_*(f,X^-)\right) \ar[d]^{I_*(f,X^-)}\\
 & H_{*}(M;\Z)
  }\] 
  }
\nd By taking the transpose of all  morphisms of chain complexes we get a similar diagram 
in cohomology made of isomorphisms:
\[\xymatrix@1{
 H^*\left(C^*(g,Y^-)\right)\ar[rd]_{I^*(g,Y^-)}&
H^*\left(C^*(f,X^-)\right) \ar[l]_(0.5){[\ga^*]} \ar[d]^{I^*(f,X^-)}\\
 & H^{*}(M;\Z)
  }\] 
  Note that a change of orientations of some unstable manifolds has the same effect on $[\ga_*]$ and on $I_*(-,-)$.
  So, the commutativity of the above diagrams is not affected.\\

 \nd 2) We can do the same for the Morse complex $C_*(f, X^+)$ which calculates the 
 relative homology. Here, we will use the stable manifolds of $X^+$ that we introduced in the beginning of Section 3. More precisely, there is a canonical process,
 similar to the one above-mentioned for the complex $C_*(f,X^-)$, which removes  
 the positive type critical points  of $f_\partial$. After removing them,
  the stable manifolds of $X^+$ associated with $crit f$ give rise to  a filtration of $M$ starting from $\partial M$:
 $$\partial M\subset M_{[1]}\subset \ldots \subset M_{[k]}\subset \ldots \subset M_{[n]}=M\,.
 $$
 Here, $M_{[k]}$ is the union of $\partial M $ and the closure of the stable manifolds
of $X^+$ converging 
 to $crit_kf$.
 The cellular homology  associated with this filtration gives a canonical isomorphism
  $$I_*(f,X^+):H_*\left(C_*(f,X^+)\right)\to H_*(M,\partial M;\Z).$$
  %is canonically isomorphic to $H_*(M,\partial M;\Z)$; denote this isomorphism by $I_*(f,X^+)$. 
  Moreover, this isomorphism is \emph{natural }
 with respect to change of function and quasi-gradient in the same sense as it is detailed in 1) above.\\
 
 \nd 3) Here comes the important point for orientations. Let $\ep_f^+$ be a choice of orientations of the stable manifolds
 of $X^+$. Since $M$ is oriented, the unstable manifolds of $X^+$ are not only co-oriented but they are also oriented.\footnote{Here, some convention has to be used, for instance: co-or(-)$\wedge$ or(-)= or($M$).}
 The latter orientations are denoted by $\ep_f^\perp$.
 
 We recall that $X^+$ is a negatively adapted quasi-gradient 
of $-f$; we denote it by $Y^-:= X^+$ when it is considered as a descending quasi-gradient of $-f$. So, we have a chain complex $C_*(-f, Y^-, \ep_{-f}^-)$ where $\ep_{-f}^-$ is determined by $\ep_f^+$ by the rule
\begin{equation}\ep_{-f}^- {=}\ep_f^\perp.
\end{equation}

By applying the functor 
$Hom(-,\Z)$ we have  its dual, a co-chain complex, $C^*(-f, Y^-,\ep_f^\perp)$. 
 By construction of $C_*^+$, we have 
 \begin{equation}\label{obvious}
\quad\eta_*: C_*(f, X^+, \ep_f^+)\mathop{\longrightarrow}^{=}C^{n-*}(-f, Y^-,\ep_f^\perp)\,.
\end{equation}
This equality means same generators and same differential; only the grading is reversed.
It induces the equality $H_*\left(C_*(f, X^+,\ep_f^+)\right)= H^{n-*}\left(C^{n-*}(-f, Y^-,\ep_f^\perp)\right)$ and by
combining it with the isomorphims  $I_{*}(f,X^+)$ and $I^{n-*}(-f,Y^-)$ we get a description at the Morse complex level of
% a natural isomorphism,
the Poincar\'e-Lefschetz isomorphism: %indeed:
$$P: H_{*}(M,\partial M;\Z)\to H^{n-*}(M;\Z)\,.$$ 
\medskip

  \begin{rien} {\bf Application to the intersection pairing.} {\rm
 We are interested in  describing 
  a pairing 
  at the
   chain level 
  $$ \si:  C_{k}(f,X^+)\otimes C_{n-k}(f,X^-)\to \Z
  $$ which induces the intersection pairing in homology.  
This is achieved in the following way. %  correcting what is said  at the end of \cite{lauden}. 
  
At the homology level the Poincar\'e-Lefschetz isomorphism  $P$ 
%which is an isomorphism $$P: H_{*}(M,\partial M;\Z)\to H^{n-*}(M;\Z),$$ 
carries the intersection product 
$$\iota:H_{*}(M,\partial M;\Z)\otimes H_{n-*}(M,\Z)\to \Z$$
to the evaluation map
$ev: H^{n-*}(M; \Z)\otimes H_{n-*}(M;\Z)\to\Z$\,.

After what was done in the previous subsection, we only have 
  to understand this evaluation map in the setting of Morse homology.
First, there is  a canonical 
evaluation map 
$$ev=<-,->: C^{n-*}(-f,Y^-)\otimes C_{n-*}(-f,Y^-)\to\Z$$ 
which on the basis elements
 is the Kronecker product. 
% - remember that $C^{n-*}_N= C^{n-*}_N(f,X^-)= Hom(C_{n-*}^N,\Z)$.
A more sophisticated way to say the same thing is to count  the  transverse
 intersection 
$W^s(x, Y^-)\cap W^u(y, Y^-)$ for every pair of critical points of the same degree, that is both in 
$crit_{k}f\cup crit^+_{k-1}f_\partial$ for some integer $k$. Here, it is essential $Y^-$ to be Morse-Smale 
for avoiding undesirable  orbits connecting points of the same degree.

We choose a generic  path\footnote{First, choose a generic path ($f_t$) in the space of functions; then, complete with
a path of quasi-gradients. For this second step, use the convexity of the set of quasi-gradients adapted to $f_t$ for a given $t$. If the function $f_t$ has a codimension-one singularity, the involved critical point needs to be a zero of $X_t$ of corank one and the previous argument still works. }
$\Gamma$ from $(f,X^-)$ to  $(-f, Y^-)$ which yields a quasi-isomorphism
$\Ga_*:C_{n-*}(f,X^-,,\ep_f^-) \to  C_{n-*}(-f, Y^-,\ep_f^\perp)$. %, well defined up to orientation signs.
Thanks to (\ref{obvious}), the desired evalution map is given by
\begin{equation}
 \si= ev\circ(\eta_*\otimes \Ga_*)\,.
 \end{equation}
If necessary, by Remark
\ref{fundrem} we may approximate $X^-$ in order to make 
mutually transverse  $W^s(x, X^+)$ and $W^u(y, X^-)$
for every $x\in crit_{k}f\cup crit^+_{k-1}f_\partial $  and $y\in crit_{n-k}f\cup 
crit^-_{n-k}f_\partial$\,.
}
\end{rien}

\nd {\sc Claim.} {\it For every pair of cycles $\al\in C_{k}(f,X^+)$
and $\beta\in C_{n-k}(f,X^-)$, the geometric formula for
$\si(\al,\beta)$ is given by counting the signed intersection 
number of the respective stable and unstable manifolds entering 
in the linear combinations forming $\al$ and $\beta$.}\\

Indeed, by  (\ref{natural}) the cycles $\beta $ and $\Ga_*(\beta)$ are homologous in $M$. Therefore,
they have the same algebraic intersection with the cycle $\al$. Notice that the frontier of the involved
invariant manifolds does not appear in this counting since it is made of invariant manifolds of less dimension.

\begin{cor} The pairing $\si$ induces the homological intersection.

\end{cor}

\end{document}